\documentclass[11pt,leqno]{amsart}
\usepackage[latin1]{inputenc}

\input{asvv.sty}

\newcommand{\ek}{{\mathrm{K}}}
\newcommand{\ee}{{\mathrm{E}}}
\newcommand{\uhp}{\mathbb{H}}
\newcommand{\RS}{{\widehat{\mathbb{C}}}} % the Riemann sphere

\renewcommand{\Im}{{\,\operatorname{Im}\,}}
\renewcommand{\Re}{{\,\operatorname{Re}\,}}
\newcommand{\arctanh}{{\,\operatorname{arctanh}\,}}
\newcommand{\D}{{\mathbb D}}

\title{Twice-punctured hyperbolic sphere with a conical singularity and
generalized elliptic integral}
\author[Anderson, Sugawa, Vamanamurthy, and Vuorinen]%
{G. D. Anderson, T. Sugawa, M. K. Vamanamurthy, and M. Vuorinen}
\date{}

\newcounter{minutes}
\divide\time by 60
\newcounter{hours}
\multiply\time by 60
\addtocounter{minutes}{-\time}
\begin{document}

\begin{center}
{\tiny \texttt{FILE:~\jobname .tex,
        printed: \number\year-\number\month-\number\day,
        \thehours.\ifnum\theminutes<10{0}\fi\theminutes}
}
\end{center}

\begin{abstract}
We describe, in terms of generalized elliptic integrals,
the hyperbolic metric of the twice-punctured sphere
with one conical singularity of prescribed order.
We also give several monotonicity properties of the metric
and a couple of applications.
\end{abstract}

\maketitle
\markboth{\textsc{ANDERSON, SUGAWA, VAMANAMURTHY, AND
VUORINEN}}{\textsc{TWICE PUNCTURED HYPERBOLIC SPHERE}}
%\markboth{\textsc{G. D. ANDERSON, T. SUGAWA, M. K. VAMANAMURTHY, AND
%M. VUORINEN}}{\textsc{POINCAR\'E DENSITY AND HYPERGEOMETRIC FUNCTIONS}}

%%%%%%%%% SECTION
%%%%%%%%% SECTION
\section{Introduction}\label{sect:intro}

The hyperbolic metric $\rho(z)|dz|$ on the thrice-punctured sphere
$\RS\setminus\{0,1,\infty\}$ is one of
the fundamental tools in complex analysis.
Indeed, for instance, the big Picard theorem can be derived
by a careful look at the metric $\rho(z)|dz|$ and the distance induced by it.
It is known that the density function $\rho(z)$ can be expressed explicitly as
$$
  \rho(z)=\frac{\pi}{8|z(1-z)|\Re\{\ek(z)\ek(1-\bar z)\}},
$$
where $\ek(z)$ is the complete elliptic integral of the first kind
given in \eqref{eq:ek} (see \cite{Agard68} or \cite{SV05}).
On the other hand, it has been recognized that generalized
elliptic integrals $\ek_a(z)$ and $\ee_a(z),$ defined in \eqref{eq:eka} and
\eqref{eq:eea} respectively,
share many properties with the original complete
elliptic integrals (cf.~\cite{AQVV00}).

In the present paper, it is shown that the hyperbolic
metric of a twice-punctured sphere with one conical singularity
of prescribed angle can be expressed in terms of these generalized complete
elliptic integrals.

\section{Hyperbolic metric with conical singularities}

A hyperbolic metric of a compact Riemann surface $R$ with conical
singularities of angle $2\pi\theta_j,~\theta_j\in[0,+\infty)\setminus\{1\},$
at points $p_j\in R,~j=1,\dots,N,$ is a conformal metric
on $R\setminus\{p_1,\dots,p_N\}$ of the form
$ds=e^{\varphi(z)}|dz|,$ where $\varphi$ is a smooth function
satisfying the Liouville equation
\begin{equation}\label{eq:Liouville}
\varDelta\varphi=4e^{2\varphi}
\end{equation}
on $R\setminus\{p_1,\dots,p_N\}$ and possessing the asymptotic behavior
\begin{equation}\label{eq:beh}
\varphi(z)=
\begin{cases}
-(1-\theta_j)\log|z-z_j|+O(1) & \quad\text{if}~\theta_j>0,\\
-\log|z-z_j|-\log(-\log|z-z_j|)+O(1) & \quad\text{if}~\theta_j=0
\end{cases}
\end{equation}
as $z\to z_j=z(p_j),$ where $z$ is a local coordinate of $R$ around $p_j.$
Note that a conical singularity of angle $0$ is called
a {\it puncture} or a {\it cusp}.

The remainder term $O(1)$ in the above is known to be continuous
at $z=z_j$ by a detailed study of the local behavior of solutions
to the Liouville equation at the isolated singularities
by Nitsche \cite{Nit57} (see also \cite{KR09}).

Heins \cite[Chap.~II]{Heins62} proved that for a compact Riemann surface $R$
of genus $g$ and finite points in it with given angles as above,
a hyperbolic metric on $R$ with the behavior described in \eqref{eq:beh}
exists {\it uniquely} as long as the condition
\begin{equation}\label{eq:GB}
2(1-g)-\sum_{j=1}^N(1-\theta_j)<0
\end{equation}
is satisfied.
This constraint comes from the Gauss-Bonnet formula.
This result was previously known by Picard \cite{Picard05} when $g=0.$
Practically, this unique metric as above is called the (complete) hyperbolic
metric of the Riemann surface $R\setminus\{p_1,\dots,p_n\}$
with conical singularities of angle $2\pi\theta_j$ at $p_j~(j=n+1,\dots, N),$
where $\theta_1=\dots=\theta_n=0<\theta_j\ne1~(j=n+1,\dots,N).$

The hyperbolic metric treated in the present paper corresponds to
the case when $R=\RS,~g=0,~N=3,~(p_1,p_2,p_3)=(0,1,\infty)$ and
$(\theta_1,\theta_2,\theta_3)=(0,0,\alpha),$ where $0\le\alpha<1.$ %2009-03-21 2\pi taken away
Note here that this case always satisfies the condition
\eqref{eq:GB}. We denote this metric by $\rho_\alpha(z)|dz|.$

When $\alpha=0,$ the metric $\rho_{0}$ is simply the
usual hyperbolic metric $\rho$ of
$\RS\setminus\{0,1,\infty\}=\mathbb{C}\setminus\{0,1\}$
(without conical singularities).
By uniqueness of the hyperbolic metric with conical singularities,
the metric admits the obvious symmetry
$\rho_\alpha(z)=\rho_\alpha(1-z)=\rho_\alpha(\bar z).$

We remark that for a M\"obius transformation $M,$
$\rho_\alpha(M(z))|M'(z)|$ gives the density
of the hyperbolic metric of $\RS\setminus\{M(0), M(1)\}$
with a conical singularity of angle $2\pi\alpha$ at $M(\infty).$
For instance, the hyperbolic metric $\tilde\rho_\alpha(z)|dz|$
of the twice-punctured sphere
$\RS\setminus\{1, \infty\}=\mathbb{C}\setminus\{1\}$
with a conical singularity of angle $2\pi\alpha$ at $0$
can be obtained by $\tilde\rho_\alpha(z)=\rho_\alpha(1/z)/|z|^2.$

%%%%%%%%% SECTION
%%%%%%%%% SECTION
\section{Generalized elliptic integrals}
The complete elliptic integrals of the first and the second kind are
defined, respectively, by
\begin{equation}\label{eq:ek}
  \ek(z)=\int_0^1\frac{dt}{\sqrt{(1-t^2)(1-zt^2)}}
\quad\text{and}\quad
  \ee(z)=\int_0^1\sqrt{\frac{1-zt^2}{1-t^2}}dt.
\end{equation}
Note that these functions can be expressed also by the hypergeometric function:
$$
  F(a,b;c;z)\!= \!{}_2 F_1(a,b;c;z)\! \equiv\!
  \sum_{n=0}^{\infty} \frac{(a,n)(b,n)}{(c,n)} \frac{z^n}{n!},\quad |z|<1,
$$
where $(a,0)=1$ and  $(a,n)= a(a+1)(a+2) \cdots (a+n-1)$ for $n\ge1, $
and $c\ne 0,-1,-2,\dots.$
Indeed,
$$
\ek(z)=\frac\pi2 F\left(\frac{1}{2}, \frac{1}{2};1;z\right)
\quad\text{and}\quad
\ee(z)=\frac\pi2 F\left(-\frac{1}{2}, \frac{1}{2};1;z\right).
$$

Let $0<a<1.$
The generalized complete elliptic integrals
of the first and the second kind with signature $1/a$
are defined, respectively, by
\begin{align}\label{eq:eka}
\ek_a(z)
&=\frac\pi2 F(a,1-a;1;z) \\
%&=\frac{\sin(\pi a)}2\int_0^1\frac{dt}{t^a(1-t)^{1-a}(1-z\,t)^a}
&=\sin(\pi a)\int_0^1\frac{t^{1-2a}dt}{(1-t^2)^{1-a}(1-z\,t^2)^a}
\notag
\end{align}
\noindent and
\begin{align}\label{eq:eea}
\ee_a(z)
&=\frac{\pi}{2} F(a-1,1-a;1;z) \\
%&=\frac{\sin(\pi a)}2\int_0^1\frac{dt}{t^a(1-t)^{1-a}\ (1-z\,t)^{a-1}}.
&=\sin(\pi a)\int_0^1\left(\frac{1-z\,t^2}{1-t^2}\right)^{1-a}t^{1-2a}dt.
\notag
\end{align}
\noindent
Here, note that $\ek_a(z)$ and $\ee_a(z)$ are defined as (single-valued)
analytic functions in $z\in\mathbb{C}\setminus[1,+\infty).$

We remark that the above definition is slightly different from the usual one.
The (traditional) generalized complete elliptic integrals
 of the first and the second kind usually
refer to $\ek_{a}(x^2)$ and $\ee_{a}(x^2)$ for $0<x<1$ in our notation.

We mean by $\ek_a'$ and $\ee_a'$ the derivatives of $\ek_a$ and $\ee_a,$
though these are often used to mean the complementary functions.
For the complementary functions, we adopt the notation
$\ek_a^*(z)=\ek_a(1-z)$ and $\ee_a^*(z)=\ee_a(1-z)$ in the present paper.

The following formula, which is a special case of Elliott's identity
(see \cite{AQVV00}),
will be used at a crucial step in the computation of the hyperbolic metric:
\begin{equation}\label{eq:Elliott}
  \ek_a^*(z)\ee_a(z)+\ee_a^*(z)\ek_a(z)-\ek_a^*(z)\ek_a(z)
  =\frac{\pi\sin(\pi a)}{4(1-a)}.
\end{equation}

Some information about
the behavior of the hypergeometric function near $z=x=1$ will be needed below.
The following result can be found in Chapter 15 of the book \cite{AS:hand}.

\begin{equation}\label{eq:viides}
\begin{cases}
F(a,b;c;1-) = \dfrac{\Gamma(c) \Gamma(c-a-b)}{\Gamma(c-a) \Gamma(c-b)},~
c>a+b,\\
F(a,b;a+b;x) = \dfrac1{B(a,b)}\log\dfrac{1}{1-x}+O(1),\\
\qquad\qquad\qquad (\text{as}~ x\to1-), \\
F(a,b;c;x) = (1-x)^{c -a -b} F(c-a,c-b;c;x),~ c<a+b.
\end{cases}
\end{equation}
Here  $B(a,b)$ denotes the beta function.

\section{Computation of $\rho_\alpha(z)$}

A relation between conformal mappings and the generalized complete elliptic
integral $\ek_a(z)$ of the first kind is given by \cite[Theorem 2.2]{AQVV00}
when the argument $z$ is real and between $0$ and $1.$
We will now give another aspect of $\ek_a(z)$ for the complex argument $z.$

As is stated in \cite[Lemma 2]{TZ03}, the hyperbolic metric of the sphere
with given conical singularities can be described in terms of solutions
to a second-order Fuchsian differential equation with regular singularities
at the cone points.
In our case, the metric is described explicitly
in terms of generalized elliptic integrals.

We begin with a general case.
Let $a$, $b$, and $c$ be real numbers.
It is a classical fact that the function
$f(z)=i\,F(a,b;a+b+1-c;1-z)/F(a,b;c;z)$ maps the upper half plane $\uhp$ onto a
curvilinear triangle bounded by three circular arcs and
having the interior angles
$(1-c)\pi$ at $f(0),~(c-a-b)\pi$ at $f(1),$ and $(b-a)\pi$ at $f(\infty),$
provided that these angles are all nonnegative and the sum is less than $\pi$
(see, for instance, \cite [pp.~206, 207]{Nehari:conf}).
Note that the segment $(0,1)$ of the real axis is mapped by $f$ to
a part of the imaginary axis
and that $f$ maps $\mathbb{C}\setminus((-\infty,0]\cup[1,+\infty))$
conformally onto the domain which is the union of $f((0,1)),$
$f(\uhp)$, and its reflection in the imaginary axis.

In the particular case when $0<a<1,~b=1-a$, and $c=1,$ the function $f$
can be written in the form $i\,\ek_a(1-z)/\ek_a(z)$, and the image $f(\uhp)$
is a circular triangle with interior angles $0$, $0$, and $|1-2a|\pi.$
More specifically, we have the following result.
\bigskip

\begin{lem}\label{lem:fa}
Let $f_a(z)=i\,\ek_a(1-z)/\ek_a(z),~0<a<1.$
Then the image $f_a(\uhp)$ of the upper half plane $\uhp$ under $f_a$
is the hyperbolic triangle $\Delta_a$
in $\uhp$ whose interior angles are $0$, $0$, and $|1-2a|\pi$
at the vertices  $0$, $\infty$, and $e^{|1-2a|\pi i/2},$ respectively.
More precisely,
$\Delta_a=\{w\in\uhp :\ 0<\Re w<\sin(\pi a),
~|2w\sin(\pi a)-1|>1\}.$
\end{lem}

\proof
%Because $f_a(z)=f_{1-a}(z)$ and the assertion is unchanged under the
%replacement $a\mapsto 1-a,$ we may assume that $0<a\le 1/2.$
Since $f_a(\uhp)$ is a Jordan domain, $f_a$ extends to a homeomorphism
from the closure of $\uhp$ onto the closure of $f_a(\uhp).$
First note that $f_a$ maps the interval $(0,1)$ onto the whole positive
imaginary axis.
Since the interior angles of $f_a(\uhp)$ at $f_a(0)$ and $f_a(1)$
are both $0,$ the boundary arcs $f_a((1,+\infty))$ and $f_a((-\infty,0))$
are contained in hyperbolic geodesics in $\uhp$ of the forms
$|w-r|=r ~(0<r)$ and $\Re w=p~(p>0),$ respectively.
In particular, the image $f_a(\uhp)$ is a hyperbolic triangle in $\uhp.$
Since we know that these two geodesics form an angle of $\theta=|1-2a|\pi,$
we find that $r(1+\cos\theta)=p$ by elementary geometry.
Thus, it is enough to show that $p=\sin(\pi a),$
which leads to the relation $r=\sin(\pi a)/(1+\cos((1-2a)\pi))=1/(2\sin(\pi a)).$

In order to make statements precise, we introduce some notation.
Let $f$ be an analytic function defined in $\C\setminus\R.$
For each $x\in\R,$ we denote by $f^{\pm}(x)$ the limit $\lim_{t\to0+}f(x\pm it)$
(if it exists).
If $f$ extends analytically to a neighborhood $V$ of $x$ as a single-valued
function on $(\C\setminus\R)\cup V,$ then we write simply
$f(x)$ as usual instead of $f^+(x)=f^-(x).$

With this notation, using \cite[(3), (27), pp.~105, 106]{Bat1},
we obtain the transformation formulas
\begin{equation}
\label{eq:-x}
\frac2\pi\ek_a(-x)=
F(a,1-a;1;-x)=(1+x)^{-a} F(a,a;1;\tfrac{x}{1+x})
\end{equation}
and
\begin{align}
\label{eq:1+x}
&\quad\ \frac2\pi\ek_a^\pm(1+x)=F^\pm(a,1-a;1;1+x)\\
&=(1+x)^{-a}\left[
\frac{\Gamma(a)}{\Gamma(2a)\Gamma(1-a)}F(a,a;2a;\tfrac{1}{1+x})
-e^{\mp \pi ai}F(a,a;1;\tfrac{x}{1+x})
\right], \notag
\end{align}
for all $x>0.$
Therefore,
\begin{equation}\label{eq:fa+}
f_a^+(-x)=i\frac{\ek_a^-(1+x)}{\ek_a^+(-x)}
=\frac{i\,\Gamma(a)}{\Gamma(2a)\Gamma(1-a)}\cdot
\frac{F(a,a;2a;\tfrac{1}{1+x})}{F(a,a;1;\tfrac{x}{1+x})}-i\,e^{\pi a\,i},
\quad x>0.
\end{equation}
Thus $\Re f_a^+(-x)=\sin(\pi a)$ for $x>0,$ as required.
\qed

\bigskip

\begin{rem}\label{rem:fa}
Since $f_a^+(x)=f_{1-a}^+(x),$ as a by-product of \eqref{eq:fa+},
we obtain the relation
\begin{align*}
&\frac{\Gamma(a)}{\Gamma(2a)\Gamma(1-a)}\cdot
\frac{F(a,a;2a;\tfrac{1}{1+x})}{F(a,a;1;\tfrac{x}{1+x})}-\cos(\pi a) \\
=~&\frac{\Gamma(1-a)}{\Gamma(2-2a)\Gamma(a)}\cdot
\frac{F(1-a,1-a;2-2a;\tfrac{1}{1+x})}{F(1-a,1-a;1;\tfrac{x}{1+x})}+\cos(\pi a)
\end{align*}
for $0<a<1$ and $x>0.$
This is equivalent to the identity
\begin{align*}
&\frac{\Gamma(a)}{\Gamma(2a)\Gamma(1-a)}\cdot F(a,a;2a;1-x)F(1-a,1-a;1;x) \\
&\ -\frac{\Gamma(1-a)}{\Gamma(2-2a)\Gamma(a)}\cdot
F(1-a,1-a;2-2a;1-x)F(a,a;1;x) \\
&\ -2\cos\pi a \cdot F(a,a;1;x)F(1-a,1-a;1;x)=0.
\end{align*}
As far as we know, this is a new identity for hypergeometric functions.
%It would be desiderable to find a proof of this identity  based only
%on the properties of the hypergeometric functions.
\end{rem}

Let $\alpha=|1-2a|.$
Since $f_a$ maps each of the intervals $(-\infty, 0)$, ~$(0,1)$,
and $(1,+\infty)$
onto a hyperbolic geodesic segment in $\uhp,$ the pull-back
$f_a^*\rho_\uhp=\rho_\uhp(f_a(z))|f_a'(z)||dz|$
of the hyperbolic (or Poincar\'e) metric $\rho_\uhp(z)|dz|=|dz|$
$/(2\Im(z))$
of the upper half plane $\uhp$,\ together with its reflection
$f_a^*\rho_\uhp(\bar z)|dz|$ defines a smooth conformal metric
on $\C\setminus\{0,1\}.$
This is the hyperbolic metric $\rho_\alpha(z)|dz|$ of
the twice-punctured sphere $\sphere\setminus\{0,1\}$
with a conical singularity of angle $2\pi\alpha$ at $\infty$
(cf.~\cite[Lemma 2]{TZ03}).
We emphasize that the curvature equation, which is equivalent to
\eqref{eq:Liouville},
\begin{equation}\label{eq:curvature}
\varDelta \log\rho_\alpha=4{\rho_\alpha}^2
\end{equation}
plays an important role in investigation of the metric.

Agard \cite{Agard68} gave a formula for $\rho_{\C\setminus\{0,1\}}=\rho_0$
in terms of complete elliptic integrals.
In the same way, we can compute $\rho_\alpha$ for $0\le\alpha<1$
with the help of the above construction.
\bigskip

\begin{thm}\label{cone}
Let $0\le\alpha<1$ and choose $0<a<1$ so that $\alpha=|1-2a|.$
The hyperbolic metric $\rho_\alpha(z)|dz|$ of the twice-punctured sphere
$\RS\setminus\{0,1\}$
with conical singularity of angle $2\pi\alpha$ at $\infty$ is given by
\begin{equation}\label{eq:rho_a}
\rho_\alpha(z)
=\frac{\pi\cos(\pi\alpha/2)}{8|z(1-z)|\Re\big(\ek_a(z)\ek_a(1-\bar z)\big)}.
\end{equation}
\end{thm}

\proof
By Gauss' contiguous relations (see (2.5.8) of \cite{AAR:special}), one obtains
$$
z(1-z)\ek_a'(z)=(1-a)\big[\ee_a(z)-(1-z)\ek_a(z)\big].
$$
Using this identity, we derive
\begin{align*}
f_a'(z)&=-i\,\frac{\ek_a'(1-z)\ek_a(z)-\ek_a(1-z)\ek_a'(z)}{(\ek_a(z))^2} \\
&=-i\,\frac{1-a}{z(1-z)}\cdot
\frac{\ee_a^*(z)\ek_a(z)+\ek_a^*(z)\ee_a(z)-\ek_a^*(z)\ek_a(z)}{(\ek_a(z))^2} \\
&=-i\,\frac{\pi\sin(\pi a)}{4z(1-z)(\ek_a(z))^2},
\end{align*}
where we have used \eqref{eq:Elliott}.
Hence, using the relation $\overline{\ek_a(z)}=\ek_a(\bar z),$
we obtain
$$
\rho_\alpha(z)=\frac{|f_a'(z)|}{2\Im f_a(z)}
=\frac{\pi\sin(\pi a)}{8|z(1-z)(\ek_a(z))^2|\Re(\ek_a(1-\bar z)/\ek_a(\bar z))},
$$
from which the required formula follows.
\qed

\bigskip
By the representation formula for $\rho_\alpha,$ we have the following.

\begin{cor}
The quantity $\rho_\alpha(z)$ is jointly continuous in $\alpha$
and $z.$
\end{cor}

\bigskip

Because the formula
\begin{equation}\label{eq:sqrt2}
\ek_a(\tfrac12)=\frac{\Gamma(\tfrac{1-a}2)\Gamma(\tfrac{a}2)\sin(\pi a)}
{4\sqrt\pi}
\end{equation}
is known (see, for instance, \cite[(4.5)]{AQVV00}), we have the
following consequence.
\bigskip

\begin{cor}
$$
\rho_\alpha(\tfrac12)
=\frac{8\pi^2}{ \left(\Gamma(\tfrac{1+\alpha}4)\right)^2%
\left(\Gamma(\tfrac{1-\alpha}4)\right)^2\cos(\tfrac{\pi\alpha}2)}.
$$
\end{cor}

The explicit formula in \eqref{eq:rho_a} of
$\rho_\alpha$ can be used to determine
the constant terms of asymptotic expansions of $\rho_\alpha$ around
singularities.

\bigskip

\begin{thm}\label{thm:asymp}
For $0<\alpha<1,$ the metric $\rho_\alpha$ satisfies
$$
\log\rho_\alpha(z)=
\begin{cases}
\log\frac1{|z|}-\log\log\frac1{|z|}-\log2+o(1)
&\quad\text{as}~ z\to 0,\\
\log\tfrac1{|z-1|}-\log\log\tfrac1{|z-1|}-\log2+o(1)
&\quad\text{as}~ z\to 1,\\
-(1+\alpha)\log|z|+
\log\frac{(\Gamma(\frac{1+\alpha}2))^2\Gamma(1-\alpha)}%
{(\Gamma(\frac{1-\alpha}2))^2\Gamma(\alpha)}+o(1)
&\quad\text{as}~ z\to \infty.
\end{cases}
$$
\end{thm}

\bigskip

\noindent
%{\sl Proof of Theorem \ref{thm:asymp}}.
\proof
Choose $a\in(0,1/2]$ so that $1-2a=\alpha.$
First we investigate $\rho_\alpha(z)$ around $z=0.$
Since the $O(1)$ term, say $w(z),$ is known to be continuous at $z=0$
(see \cite[Satz 1]{Nit57} or \cite[Theorem 1.1]{KR09}),
it suffices to show that $w(0)=\log 2.$
By \eqref{eq:viides}, for $x>0$ we have
$$
\ek_a(x)=\frac\pi2+O(x)\quad\text{and}\quad
\ek_a(1-x)=\frac{\sin(\pi a)}{2}\log\frac1{x}+O(1)
$$
as $x\to 0+.$
Substitution of these formulas into \eqref{eq:rho_a} yields
$w(0)=\log2$ as required.
The corresponding result for $z=1$ follows from the previous one
by the symmetry $\rho_\alpha(1-z)=\rho_\alpha(z).$

Finally, we consider the case $z\to\infty.$
By the general property of conical singularities, one has the expression
$\log\rho_\alpha(z)=(2a-2)\log|z|+v(z),$ where $v(z)$ is a continuous function
near $z=\infty$ (see \cite[Satz 1]{Nit57} or \cite[Theorem 1.1]{KR09}).
For $x>0,$ by \eqref{eq:-x}, \eqref{eq:1+x}, and \eqref{eq:viides}, we have
\begin{equation}\label{eq:kax}
\frac2\pi\ek_a(-x)=\frac{\Gamma(1-2a)}{(\Gamma(1-a))^2}x^{-a}(1+o(1))
\end{equation}
and
\begin{align}\label{eq:ka1x}
&\quad\
\frac2\pi\Re \ek_a^\pm(1+x) \\
\notag
&=(1+x)^{-a}\left[\frac{\Gamma(a)}{\Gamma(2a)\Gamma(1-a)}
F(a,a;2a;\tfrac1{1+x})-\cos(\pi a)\,F(a,a;1;\tfrac x{1+x})
\right] \\
\notag
&=\left[\frac{\Gamma(a)}{\Gamma(2a)\Gamma(1-a)}
-\cos(\pi a)\frac{\Gamma(1-2a)}{(\Gamma(1-a))^2}\right]\,x^{-a}(1+o(1)) \\
\notag
&=\frac{\Gamma(a)}{2\Gamma(2a)\Gamma(1-a)}\,x^{-a}(1+o(1)) \\
\notag
&=\frac{(\Gamma(a))^2\sin(\pi a)}{2\pi \Gamma(2a)}\,x^{-a}(1+o(1))
\notag
\end{align}
as $x\to+\infty.$
Combining \eqref{eq:kax} and \eqref{eq:ka1x} with \eqref{eq:rho_a}, we see that
$$
\rho_\alpha(-x)=\frac{(\Gamma(1-a))^2\Gamma(2a)}{(\Gamma(a))^2\Gamma(1-2a)}\,
x^{2a-2}(1+o(1)),\quad x\to+\infty,
$$
which implies that
$v(\infty)=\log(\Gamma(1-a))^2\Gamma(2a)/\Gamma(a)^2\Gamma(1-2a)).$
This is equal to the required constant term.
\qed

\bigskip

Lehto, Virtanen and V\"ais\"al\"a \cite{LVV59} proved the useful inequality
$\rho_{0}(-|z|)\le \rho_{0}(z)$ for all $z\in\C\setminus\{0,1\}.$
Later on, Weitsman \cite{Weit79} proved a monotonicity property
of the hyperbolic metric on a circularly symmetric domain, which
means that $\rho_{0}(r\,e^{i\theta})$ is a non-increasing function
of $\theta$ in $0<\theta<\pi$ for a fixed $r>0$ for the particular domain
$\C\setminus\{0,1\}.$
We can deduce the same result for $\rho_\alpha$
by employing the method developed in \cite{LVV59}.
\bigskip

\begin{thm}\label{thm:LVV}
For $0\le\alpha<1$ and fixed $r>0,~\rho_\alpha(r\,e^{i\theta})$ is
a non-increasing (non-decreasing) function of $\theta$ in $0<\theta<\pi~
(-\pi<\theta<0).$
In particular, the inequalities
$\rho_\alpha(-|z|)\le\rho_\alpha(z)\le\rho_\alpha(|z|)$
hold for each $z\in\C\setminus\{0,1\}.$
\end{thm}

\proof
It is enough to show the assertion by assuming that $0<a<\frac 12.$
(The case $a=\frac 12$ can be treated similarly with the special relation
$\rho_{0}(1/z)=\rho_{0}(z)|z|^2$ being taken into account.)
By the obvious symmetry $\rho_\alpha(\bar z)=\rho_\alpha(z),$ it is enough to
prove the inequality
$\rho_\alpha(r\,e^{i\theta_1})\ge\rho_\alpha(r\,e^{i\theta_2})$
for $0\le\theta_1<\theta_2\le\pi.$
Let $\lambda_1(z)=\rho_\alpha(e^{-i\theta_0}z)$
and $\lambda_2(z)=\rho_\alpha(e^{i\theta_0}z),$
where $\theta_0=(\theta_1+\theta_2)/2.$
Consider now the function $h(z)=\log\lambda_1(z)-\log\lambda_2(z).$
Then $h$ is smooth in $\C\setminus\{0,e^{i\theta_0},e^{-i\theta_0}\}$ and,
by the above symmetry, $h=0$ on $\R\setminus\{0\}.$

We will show that $h(z)\ge 0$ for $z\in\uhp.$\
To this end, we first observe the asymptotic behavior of $h(z).$
It is easy to see that $h(z)\to+\infty$ as $z\to e^{i\theta_0}.$
By Theorem \ref{thm:asymp}, we also have
$h(z)\to0$ as $z\to\infty$ or $z\to0.$
Therefore, the set $\{z\in\uhp\setminus\{e^{i\theta_0}\}: h(z)\le -\varepsilon\}$
is compact for each $\varepsilon>0.$
Suppose now that $h<0$ somewhere in $\uhp.$
Then, there would be a
minimum point $z_0$ for $h$ in $\uhp\setminus\{e^{i\theta_0}\}.$
Then $\varDelta h(z_0)\ge0$ by minimality.
On the other hand, the inequality $h(z_0)<0$ would imply
$\lambda_1(z_0)<\lambda_2(z_0).$
Hence, by \eqref{eq:curvature},
$$
\varDelta h(z_0)=4\lambda_1(z_0)^2-4\lambda_2(z_0)^2<0,
$$
which would be impossible.
Thus, we have shown that $h(z)\ge0$ for $z\in\uhp.$
We now take the point $z_0=r\,e^{i(\theta_2-\theta_1)/2}.$
Then $0\le h(z_0)=\log\rho_\alpha(r\,e^{-i\theta_1})-\log\rho_\alpha(r\,e^{i\theta_2})
=\log\rho_\alpha(r\,e^{i\theta_1})-\log\rho_\alpha(r\,e^{i\theta_2}),$
and thus, $\rho_\alpha(r\,e^{i\theta_1})\ge\rho_\alpha(r\,e^{i\theta_2}).$
\qed

\medskip

The hyperbolic distance on $\sphere\setminus\{0,1\}$
induced by $\rho_\alpha$ is defined, as usual, by
$$
d_\alpha(z_1,z_2)=\inf_\gamma\int_\gamma\rho_\alpha(z)|dz|,
$$
where $\gamma$ runs over all the rectifiable paths $\gamma$ connecting
$z_1$ and $z_2$ in $\RS\setminus\{0,1\}.$

As a corollary of Theorem \ref{thm:LVV}, we derive a lower estimate for
the hyperbolic distance.
\bigskip

\begin{cor}
For $0<a<1$ and $z_1,~z_2\in\RS\setminus\{0,1\}$ with $|z_1|\le|z_2|,$
the following inequality holds:
\begin{equation}\label{eq:lower}
d_\alpha(z_1,z_2)\ge d_\alpha(-|z_1|,-|z_2|)
=\int_{|z_1|}^{|z_2|}\rho_\alpha(-t)dt.
\end{equation}
\end{cor}

\bigskip

We can compute the last integral by the following result.
\bigskip

\begin{thm}\label{hypdist}
Let $\alpha=|1-2a|$ for $0<a<1.$
The formula
$$
\int_x^y\rho_\alpha(-t)dt=\Phi_a(y)-\Phi_a(x)
$$
holds for $0<x<y,$ where
$$
\Phi_a(x)=-\frac12\log\left(
\frac{\Gamma(a)}{\Gamma(2a)\Gamma(1-a)}
\frac{F(a,a;2a;\tfrac{1}{1+x})}{F(a,a;1;\tfrac{x}{1+x})}-\cos(\pi a)
\right).
$$
\end{thm}

\proof
One can proceed almost as in the proof of \cite[Lemma 5.1]{SV05}.
We can write $f_a^+(-x)$ in the form $i\,u(x)+\sin(\pi a)$ for $x>0$
by \eqref{eq:fa+}.
Since $\rho_\alpha(-t)=|(f_a^+)'(-t)|/2\Im f_a^+(-t)=-u'(t)/2u(t),$
we obtain
$$
\int_x^y\rho_\alpha(-t)dt
=-\int_x^y\frac{u'(t)}{2u(t)}dt=\frac12\log\frac{u(x)}{u(y)}
=\Phi_a(y)-\Phi_a(x).
$$
\qed

\bigskip

Note that when $a\ne \frac12,$
$$
\Phi_a(\infty)=-\frac12\log\cos(\pi a)
$$
is positive and finite, whereas $\Phi_{1/2}(\infty )=\infty. $
\bigskip

\begin{rem}
More generally, the $\rho_\alpha$-distance between $z_1$ and $z_2$
in $\overline{\uhp}\setminus\{0,1,\infty\}$ can be expressed by
$$
d_\alpha(z_1,z_2)=\arctanh\left|
\frac{f_a(z_2)-f_a(z_1)}{f_a(z_2)-\overline{f_a(z_1)}}
\right|,
$$
where $0<a<1$ is chosen so that $\alpha=|1-2a|$ and $f_a$ is given in
Lemma \ref{lem:fa}.
Indeed, by construction, $f_a$ is an isometric embedding
of $(\uhp,\rho_\alpha)$ into $(\uhp,\rho_\uhp)$
and its image $\Delta_a$ is (hyperbolically) convex in $\uhp.$
Therefore, the geodesic segment joining $z_1$ and $z_2$ in $\RS\setminus\{0,1\}$
with respect to $\rho_\alpha$ is contained in the closure of $\uhp$
and its image under $f_a$ is the hyperbolic geodesic joining
$f_a(z_1)$ and $f_a(z_2).$
It is well known that the hyperbolic distance between two points $w_1$
and $w_2$ in $\uhp$ is given by $\arctanh|(w_2-w_1)/(w_2-\overline{w_1})|,$
and the above formula follows.
\end{rem}

Finally, we mention monotonicity of $\rho_\alpha(z)$ with respect to the
parameter $\alpha.$
\bigskip

\begin{prop}
The density $\rho_\alpha(z)$ is non-increasing in $0\le\alpha<1$ for a fixed
$z\in\C\setminus\{0,1\}.$
\end{prop}

Though this result is contained in \cite[Prop.~2.4]{ST05}
as a special case, we give a proof for convenience of the reader.
The assertion is established
by a simple application of the Schwarz-Pick-Ahlfors lemma
(cf. \cite{Ahlfors:conf}).
Here, we employ the same technique as in Theorem \ref{thm:LVV}.

\proof
For a given pair $\alpha, \alpha'$ with $0\le\alpha<\alpha'<1,$
we consider the function
$h=\log\rho_{\alpha'}-\log\rho_\alpha$ in $\C\setminus\{0,1\}.$
By Theorem \ref{thm:asymp}, the function $h$ extends continuously to $0$ and $1$
if we set $h(0)=h(1)=0$, and has the asymptotic behavior
$h(z)=(\alpha-\alpha'+o(1))\log|z|$ as $z\to\infty.$
Therefore, if $h$ takes a positive value, there is a point
$z_0\in\C\setminus\{0,1\}$ at which $h$ attains its (positive) maximum.
Then $\varDelta h(z_0)\le0.$
On the other hand, by \eqref{eq:curvature},
$$
\varDelta h(z_0)=4\rho_{\alpha'}(z_0)^2-4\rho_\alpha(z_0)^2>0,
$$
which is a contradiction.
Hence, we conclude that $h(z)\le0,$ in other words,
$\rho_\alpha(z)\ge\rho_{\alpha'}(z)$ for $z\in\C\setminus\{0,1\}.$
\qed

\bigskip

\begin{rem}
The expression $\rho(a,z)\equiv \rho_{|1-2a|}(z)$ as in \eqref{eq:rho_a}
can be viewed as a smooth function in $(a,z)\in(0,1)\times(\C\setminus\{0,1\}).$
Then it has the obvious symmetry $\rho(1-a,z)=\rho(a,z).$
By the above theorem, $\rho(a,z)$ attains its maximum at $a=\frac 12$ for
a fixed $z.$
In particular, by this observation we obtain
$(\partial\rho/\partial a)(\frac 12,z)=0$.
We also see that $\rho_\alpha(z)\to0$ as $\alpha\to1$ from \eqref{eq:rho_a}.
This corresponds to the well-known fact that the twice-punctured sphere
$\RS\setminus\{0,1\}$ does not carry a hyperbolic metric.
\end{rem}
\bigskip

\section{Applications}

We conclude the present note with a few applications of our metric
$\rho_\alpha.$
Since no concrete estimates for $\rho_\alpha$ are given so far,
we will give only general principles to refine classical results.

If a meromorphic function $f$ on the unit disk $\D$ does not assume the
three points $0,1$ and $\infty,$ then the principle of hyperbolic metric
gives us the inequality $f^*\rho\le \rho_\D,$ namely;
$$
\rho(f(z))|f'(z)|\le\frac1{1-|z|^2},\quad z\in\D.
$$
The classical theorems of Picard and Schottky follow essentially from the above
inequality (see, for example, \cite[\S 1-9]{Ahlfors:conf}).
We can now relax the assumption about the omitted values as in the following.

\begin{thm}
Let $f$ be a meromorphic function on the unit disk omitting the two values
$0$ and $1.$
Suppose that every pole of $f$ is of order at least $k\ge2.$
Then the following inequality holds:
$$
\rho_{1/k}(f(z))|f'(z)|\le \frac1{1-|z|^2},\quad z\in\D.
$$
\end{thm}

\proof
Let $\alpha>1/k.$
We may assume that $f$ is not constant.
Let $\lambda$ be the pull-back metric $f^*\rho_\alpha$ of $\rho_\alpha$
under $f.$
Then it is easily verified that the Gaussian curvature of $\lambda$
is $-4$ off the set of poles and branch points of $f.$
Let $z_0$ be a pole of $f.$
Then the order $m$ of the pole at $z_0$ is at least $k$ by assumption.
In view of \eqref{eq:viides}, we have
\begin{align*}
\log\lambda(z)&=-(1+\alpha)\log|f(z)|+\log|f'(z)|+O(1) \\
&=\big[ (1+\alpha)m-(m+1)\big]\log|z-z_0|+O(1) \\
&=\left(\alpha m-1\right)\log|z-z_0|+O(1)
\end{align*}
as $z\to z_0.$
Since $\alpha m>1,$ we see that $\lambda(z_0)=0.$
Thus $\lambda$ is an ultrahyperbolic metric on $\D$ in the sense of
Ahlfors \cite{Ahlfors:conf}.
Thus, Ahlfors' lemma now yields $f^*\rho_\alpha\le\rho_\D.$
Taking the limit as $\alpha\to1/k,$ we obtain the required inequality.
\qed

Knowledge about the hyperbolic metric $\rho=\rho_0$ of the thrice-punctured
sphere $\C\setminus\{0,1\}$ has led to various useful estimates for
the hyperbolic metric of a general plane domain
(see, for instance, \cite{BP78} or \cite{SV05}).
We now use $\rho_\alpha$ instead of $\rho_0$ to obtain similar
estimates for the hyperbolic metric with conical singularities.

\begin{thm}
Let $\Omega$ be a subdomain of the Riemann sphere $\RS$ with
$\infty\in\Omega$ such that $\RS\setminus\Omega$ contains at least two points.
Let $\lambda$ be a conformal metric on $\Omega$ with conical singularities
of angle less than $2\pi.$
Suppose that $\lambda$ has a conical singularity of angle $2\pi\alpha>0$
and that for each $w_0\in\partial\Omega,$
$|z-w_0|\log(1/|z-w_0|)\lambda(z)$ is bounded away from $0$
in $V\cap\Omega$ for a neighborhood $V$ of $w_0$
if $w_0$ is isolated in $\partial\Omega$
and $|z-w_0|\log(1/|z-w_0|)\lambda(z)\to+\infty$ as $z\to w_0$ in $\Omega$
otherwise.
Then
$$
\lambda(z)\ge\sup_{w_0,w_1\in\partial\Omega}
\frac1{|w_1-w_0|}\rho_\alpha\left(\frac{z-w_0}{w_1-w_0}\right),
\quad z\in\Omega\setminus\{\infty\}.
$$
\end{thm}

\proof
We follow the argument used by Heins \cite[\S 20]{Heins62}.
First note that the set $S$ of conical singularities of $\lambda$
can be characterized as
$\{z\in\Omega\setminus\{\infty\}: \lambda(z)=\infty\}\cup\{\infty\}.$
Pick $\alpha'\in(\alpha,1)$ and
fix a pair of distinct points $w_0$ and $w_1$ in $\partial\Omega.$
Set $\mu(z)=\rho_{\alpha'}((z-w_0)/(w_1-w_0))/|w_1-w_0|$ and let
$$
v=\max\{\log\mu-\log\lambda,~0\}.
$$
Then $v$ is subharmonic on $\Omega\setminus S$ and vanishes
in a neighborhood of $S.$
Moreover, $v=0$ near every boundary point of $\Omega$ except possibly
for $w_0$ and $w_1.$
If $w_j$ is not isolated, then this is still valid.
Otherwise, by the local behavior of solutions to the Liouville equation
around an isolated singularity due to Nitsche \cite{Nit57},
we see that $v$ can be extended continuously to the point $w_j.$
Recall now the following fact:
{\it
Suppose that $u$ is a continuous function on an open neighborhood $D$ of
a point $a$ and subharmonic on $D\setminus\{a\}.$
Then $u$ is subharmonic on $D.$}
Thus $v$ is subharmonic on $\Omega'=\Omega\cup\{w_j: w_j ~\text{is isolated}\}$
and vanishes in a neighborhood of $\partial\Omega'\cup S.$
We now appeal to the maximum principle to conclude that $v=0$ in $\Omega,$
which means $\mu\le\lambda.$
The proof is now complete.
\qed
\bigskip

%The following are problems for possible further contents of the present paper.
%Your suggestions are most appreciated.
%(This part will be removed in the final form of the paper.)

%\begin{open}

%(1) Find the value $\rho_\alpha(-1).$

%(2) Is the identity in Remark \ref{rem:fa} previously known?
%\end{open}

%\bibliographystyle{amsplain}
%\bibliography{papers}

\def\cprime{$'$} \def\cprime{$'$} \def\cprime{$'$}
\providecommand{\bysame}{\leavevmode\hbox to3em{\hrulefill}\thinspace}
\providecommand{\MR}{\relax\ifhmode\unskip\space\fi MR }
\providecommand{\MRhref}[2]{%
  \href{http://www.ams.org/mathscinet-getitem?mr=#1}{#2}
}
\providecommand{\href}[2]{#2}

\bigskip

%\vspace{1in}

\noindent
ANDERSON: \\
  Department of Mathematics \\
Michigan State University \\
     East Lansing, MI 48824, USA \\
      email: {\tt anderson@math.msu.edu}\\
     FAX: +1-517-432-1562\\ [1mm]

\noindent
SUGAWA:\\
     Graduate School of Information Sciences \\
     Tohoku University \\
     Aoba-ku, Sendai, 980-8579 JAPAN\\
     email: {\tt sugawa@math.is.tohoku.ac.jp}\\
     FAX: +81-22-795-4654\\ [1mm]

\noindent
VAMANAMURTHY:\\
   Department of Mathematics \\
    University of Auckland \\
    Auckland, NEW ZEALAND\\
     email: {\tt vamanamu@math.auckland.nz}\\
FAX: +649-373-7457\\ [2mm]

\noindent
VUORINEN:\\
     Department of Mathematics \\
     University of Turku \\
     Vesilinnantie 5\\
     FIN-20014, FINLAND\\
     e-mail: ~~{\tt vuorinen@utu.fi}\\
     FAX: +358-2-3336595\\

\end{document}